# New Type of Soliton Equation Described Some Statistical Distributions and Nonlinear Equations Unified Quantum Statistics

Yi-Fang Chang

Department of Physics, Yunnan University, Kunming, 650091, China

(e-mail: yifangchang1030@hotmail.com)

Abstract: We proposed a new type of soliton equation, whose solutions may describe some statistical distributions, for example, Cauchy distribution, normal distribution and student's t distribution, etc. The equation possesses two characters. Further, from an extension of this type of equation we may obtain the exponential distribution, and the Fermi-Dirac distribution in quantum statistics. Moreover, by using the method of the soliton-solution, the nonlinear Klein-Gordon equation and nonlinear Dirac equations may derive Bose-Einstein and Fermi-Dirac distributions, respectively, and both distributions may be unified by the nonlinear equation.

**1.Introduction**

In the nonlinear theory, the soliton is already applied to many regions [1,2]. Some relations among the nonlinear equations and their soliton and distributions are discussed widely. But these discussions fasten on the few known nonlinear Schr?dinger equation [3-13], the sine-Gordon equation [14-16] and the Korteweg–de Vries equation [7,9], etc.

Recently, Ablowitz, et al., found that general initial pulses are essentially solitons of the classical nonlinear Schr?dinger (NLS) equation [17]. Zhong and Belić, et al., constructed exact periodic wave solutions to the generalized two-dimensional nonlinear Schr?dinger equation with distributed dispersion, nonlinearity, and gain coefficients [18], and obtained exact spatiotemporal periodic traveling wave solutions to the generalized (3+1)-dimensional nonlinear Schr?dinger equation with distributed coefficients [19]. Derevyanko, et al., researched the random input problem for a nonlinear system modeled by the integrable one-dimensional self-focusing nonlinear Schr?dinger equation, and obtained two models and the distribution of minimal pulse width required for a soliton generation [20]. Pankratov investigated the noise self-pumping effect in long Josephson junctions (fluctuating solitons) modeled by a sine-Gordon equation [21].

**2.New type of soliton equation described some statistical distributions**

We propose a type of soliton equation [22]:

$$\varphi(\varphi_{xx} - \varphi_{tt}) - a(\varphi_x^2 - \varphi_t^2) + b\varphi^{n+1} = 0. \qquad (1)$$

By the method of the soliton-solution, let $\eta = \dfrac{x - ut}{\sqrt{1 - u^2}}$, Eq.(1) becomes an ordinary differential equation:

$$\varphi\varphi''-a(\varphi')^2+b\varphi^{n+1}=0. \tag{2}$$

$$\frac{d\varphi}{d\eta}=\varphi'=\varphi^a[-\int 2b\varphi^{n-2a}d\varphi+C_1]^{1/2}. \tag{3}$$

When $n \neq 2a-1$,

$$\varphi'=\varphi^a[C_1-\frac{2b\varphi^{n-2a+1}}{n-2a+1}]^{1/2}; \tag{4}$$

When $n = 2a-1$,

$$\varphi'=\varphi^a[C_1-2b\ln|\varphi|]^{1/2}. \tag{5}$$

(A). If $a=2$, $n=2$, $\varphi'=\varphi[C_1\varphi^2+2b\varphi]^{1/2}$,

$$\varphi=\frac{2b}{b^2(\eta+C_2)^2-C_1}. \tag{6}$$

$\varphi \to 0$ for $\eta \to \pm\infty$. It is a soliton solution. When $C_2=-\mu$, $b=\frac{2\pi}{\lambda}>0$ and $C_1=-4\pi^2$, if $\eta$ is a continuous random variable, $\varphi=\frac{\lambda}{\pi}\frac{1}{(\eta-\mu)^2+\lambda^2}$ will be Cauchy distribution density [23].

(B). If $a=1$, $n=1$,

$$\varphi=\exp[\frac{C_1-b^2(\eta+C_2)^2}{2b}]. \tag{7}$$

It is also a soliton solution. When $C_2=-\mu$, $b=\frac{1}{\sigma^2}$ and $C_1=\frac{1}{\sigma^2}\ln(\frac{1}{\sigma\sqrt{2\pi}})$,

$\varphi=\frac{1}{\sigma\sqrt{2\pi}}\exp[-\frac{(\eta-\mu)^2}{2\sigma^2}]$ is a normal distribution density [23].

(C). If $a=n=\frac{v+3}{v+1}$, $\varphi'=\varphi^{(v+3)/(v+1)}[C_1+b(v+1)\varphi^{-2/(v+1)}]^{1/2}$,

$$\varphi=[\frac{b^2(\eta+C_2)^2-C_1}{b(v+1)}]^{-(v+1)/2}. \tag{8}$$

It is still a soliton solution. When $C_2=0$, $b=\frac{v+1}{v}D^{-2/(v+1)}$, $C_1=-\frac{(1+v)^2}{v}D^{-4/(v+1)}$ and

$D=\frac{1}{\sqrt{v\pi}}\frac{\Gamma[(v+1)/2]}{\Gamma(v/2)}$, $\varphi=\frac{1}{\sqrt{v\pi}}\frac{\Gamma[(v+1)/2]}{\Gamma(v/2)}(1+\frac{\eta^2}{v})^{-(v+1)/2}$ is a student's t distribution.

All of these distributions are bell form, which are also soliton form.

When b=0 and a=1, Eq.(2) becomes $\varphi\varphi''-(\varphi')^2 = 0$, and $\varphi = C_0 e^{-c\eta}$ is the exponential distribution. But, this is not the soliton form.

Eq.(1) possesses following character: 1). If $\varphi \to y^p$ (p is constant), (1) becomes:

$$y(y_{xx} - y_{tt}) - (ap - p + 1)(y_x^2 - y_t^2) + (b/p)y^{p(n-1)+2} = 0. \tag{9}$$

The form of Eq.(1) is the same, and corresponding solitons and distributions are also invariant. Only for $p = (1-a)^{-1}$ Eq.(1) is simplified to

$$(y_{xx} - y_{tt}) + b(1-a)y^{(n-a)/(1-a)} = 0. \tag{10}$$

2). It possesses Lorentz invariance, and is relativistic covariant.

### 3. Extension of this type of equation

Further, Eq.(1) may be extended to more general nonlinear equation:

$$\varphi(\varphi_{xx} - \varphi_{tt}) - a(\varphi_x^2 - \varphi_t^2) = F(\varphi). \tag{11}$$

Eq.(11) becomes the ordinary differential equation:

$$\varphi\varphi'' - a(\varphi')^2 = F(\varphi). \tag{12}$$

$$\frac{d\varphi}{d\eta} = \varphi^a [\int 2F(\varphi)\varphi^{1-2a} d\varphi + C_1]^{1/2}. \tag{13}$$

When $F(\varphi)$ is various forms, we can derive many ordinary differential equations, in which some may be solved by the elliptic function. If $F(\varphi)$ is polynomial, for example, when a=1 and $F(\varphi) = f\varphi^4 - (b/2)\varphi^3$, for $C_1 = b^2/4f > 0$,

$$\varphi = \frac{b^2/2f}{\exp[-\sqrt{b^2/2f}(\eta + C_2)] + b}. \tag{14}$$

For b=1, this is Fermi-Dirac distribution in quantum statistics.

For the special cases a=0, Eq.(11) is similar with Eq.(10):

(a). When $F(\varphi) = \varphi \sin \varphi$, it is sine-Gordon equation.

(b). Eq.(11) can include the relativistic extension $\varphi_{xx} - \varphi_{tt} + v|\varphi|^2 \varphi = 0$ of cubic Schrodinger equation, and Higgs equation $\varphi_{xx} - \varphi_{tt} + m_0^2 \varphi - f^2 \varphi^3 = 0$.

(c). When $F(\varphi)$ is polynomial, Eq.(12) may be simplified equations of KdV equation and cubic Schrodinger equation, etc., $\varphi\varphi'' = a\varphi^2 - b\varphi^n$, here n=3 or 4.

(d). When $F(\varphi) = f\varphi^2 - b\varphi^4$, for $C_1 = 0$, f>0 and b>0,

$$\varphi = \sqrt{2f/b} \sec h(\sqrt{f}\eta + C). \tag{15}$$

For $C_1 = -f^2/2b$, f<0 and b<0,

$$\varphi = \sqrt{f/b} th(\sqrt{-b^2/2f}\eta + C). \tag{16}$$

Both are respectively the simplest soliton solutions of different forms.

**4. Nonlinear equation unified quantum statistics**

When $a=0$ and $F(\varphi) = m^2\varphi^2 + a\varphi^4$ for Eq.(11), it is namely the nonlinear Klein-Gordon equation:

$$(\varphi_{xx} - \varphi_{tt}) - m^2\varphi = a\varphi^3. \tag{17}$$

This may derive the kink. If Eq.(17) is extended to

$$(\varphi_{xx} - \varphi_{tt}) - m^2\varphi = a\varphi^3 + b^2\varphi^5. \tag{18}$$

Let $\eta = \dfrac{x - ut}{\sqrt{1-u^2}}$, Eq.(18) becomes an ordinary differential equation, and integral is made,

$$\frac{d\varphi}{d\eta} = (m^2\varphi^2 + \frac{1}{2}a\varphi^4 + \frac{1}{3}b^2\varphi^6 + C_0)^{1/2}. \tag{19}$$

Let $C_0 = 0$ and $b = \sqrt{3}a/4m$, then

$$|\varphi|^2 = \frac{2m^2}{e^{-2m\eta + C} - (a/2)}. \tag{20}$$

If $a=2$, Eq.(20) is namely Bose-Einstein (BE) distribution. The square of the wave function is just the probability [24].

Based on the Heisenberg's unified equation [25]

$$\gamma_\mu \partial_\mu \psi - l_0^2 \psi(\psi^+\psi) = 0, \tag{21}$$

which possesses the rest mass, it is namely a nonlinear Dirac equation

$$\gamma_\mu \partial_\mu \psi + m\psi - b\psi(\psi^+\psi) = 0. \tag{22}$$

In quantum mechanics, the probability density is $\rho = \psi^+\psi = \bar{\psi}\gamma_4\psi$, and $\psi^+\psi = 1 - \psi\psi^+$.

By a similar method of the soliton-solution, let $\eta' = \dfrac{\gamma_\alpha x_\alpha - u\gamma_0 t}{1+u}$, Eq.(22) becomes an ordinary differential equation

$$\frac{d\psi}{d\eta'} = b\psi(\psi^+\psi - \frac{m}{b}).$$  (23)

Its soliton solution is

$$\rho = \psi^+\psi = \frac{m/2b}{e^{2m\eta'+C}+1}.$$  (24)

It is just Fermi-Dirac (FD) distribution. Therefore, the quantum statistics (20) and (24) of bosons and fermions are obtained from the corresponding nonlinear dynamical equations (18) and (22).

Further, we extend Eq.(22) to

$$\gamma_\mu\partial_\mu\psi + m\psi = nb\psi(\psi^+\psi),$$  (25)

so $$\rho = \psi^+\psi = \frac{m/2b}{e^{2m\eta'+C}+n}.$$  (26)

Let n=1, -1 or 0, i.e., the nonlinear terms are a reversed sign or zero, Eq.(26) are FD, BE or Maxwell-Boltzmann (MB) distributions, respectively.

This is not only a supersymmetry in which fermions and bosons are complete symmetrical, but also agree with the experiments where the multiplicity and its distribution, etc., are independent of the types of particles. We proposed that the extensive chaos theory may describe the multiparticle production and the extensive air showers at high energy [26,24,27]. In the other aspect, based on the experiments and theories at high energy, we discussed a possibility: FD and BE statistics will be unified at high energy [28,24] and some possible tests of the inapplicability of Pauli's exclusion principle [29-32].

The nonlinear equations may obtain many statistical distributions. This can contact statistical distributions, differential equations and physical dynamics.